\documentstyle[A4,12pt]{article}
\newcommand \backs {\hspace{-1.5mm}}
\newcommand \goback {\hspace{-2.5 mm}}
\newcommand \halfsp {\hspace{0.5 mm}}
\newcommand \qed {\rule{2mm}{2mm}}
\newcommand \car {\hspace{-0.5 mm} \rule[-1.3mm]{0.4mm}{4.8 mm} 
                                  \hspace{-0.5 mm}}
\newcommand \carr {\hspace{-0.5 mm} \rule[-1.3mm]{0.4mm}{4.8 mm} 
		\hspace{-0.9 mm}}
\newcommand \carl {\hspace{-1 mm} \rule[-1.3mm]{0.4mm}{4.8 mm} 
		\hspace{-0.5 mm}}
\newcommand \balk {\rule{3mm}{0.3mm}}
\newcommand \ls {\hspace {0.1 cm}}
\newcommand \sq {\hspace {0.3 cm}}

\newcommand \FiniteA {{\em {fin}}(A)}
\newcommand \FiniteB {{\em {fin}}(B)}
\newcommand \FiniteS {{\em {fin}}(S)}
\newcommand \final {{\em {fin}}(\alpha)}
\newcommand \finga {{\em {fin}}(\gamma)}
\newcommand \PowerS {{\cal{P}}(S)}
\newcommand \oneone {{\hspace{-0.3mm}}{\rule[0.8mm]{1mm}{0.1mm}}
                                       {\hspace{-0.3mm}}}
\newcommand \onone {{\hspace{-0.3mm}}{\rule[1.2mm]{1.2mm}{0.1mm}}
                                     {\hspace{-1.3mm}}}
\newcommand \inj {\hbox{1{\onone}1 }}
\newcommand \Seq {{\em {seq}}^{{\em{1}}{\oneone}{\em{1}}}}
\newcommand \ief {${}^{{\em{1}}{\oneone}{\em{1}}}\ $}
\newcommand \so {s_{0}}
\newcommand \se {s_{1}}
\newcommand \sn {s_{n}}
\newcommand \si {s_{i}}
\newcommand \sa {s_{\alpha}}
\newcommand \sbe {s_{\beta}}
\newcommand \smu {s_{\mu}}
\newcommand \sio {s_{\iota}}
\newcommand \cirk {{\hspace{0.8mm}}\put(0,3){\circle{2}}
		{\hspace{1mm}}}
\begin{document}
\begin{center}
{\large{\bf  Consequences of Arithmetic for }}\\
{\large{\bf Set Theory}}
\end{center}
\vspace {1.1  cm}
\begin{center}
\small Lorenz HALBEISEN \footnote{Parts of this work are of the first 
	author's Diplomarbeit 
at the ETH Z\"urich.
He is grateful to his supervisor, Professor H.\halfsp L\"auchli.}\\
Department of Mathematics, \\ ETH Z\"urich, Switzerland
\end{center}
\vspace {0.1 cm}
\begin{center}
\small Saharon SHELAH \footnote{Research partially supported 
by the Basic Research Fund, Israeli Academy;
Publ.\halfsp No. 488}\\ 
Institute of Mathematics, \\ Hebrew University 
	Jerusalem, Israel
\end{center}
%
%
%
%
\vspace {0.6 cm}

\begin{abstract}
In this paper, we consider certain cardinals in
ZF (set theory without AC, the Axiom of Choice). In ZFC (set theory with 
AC), given  any cardinals $\cal{C}$ and $\cal{D},$ either  
$\cal{C} \leq \cal{D}$ or
$\cal{D} \leq \cal{C}.$ 
However, in ZF this is no longer so. 
For a given infinite set $A$ consider  $\Seq (A)$,
the set of all sequences of $A$ without repetition. 
We compare \carr$\Seq (A)$\carl, the cardinality of this set,
to \car${\cal {P}}(A)$\carl, the cardinality of the power set 
of $A$. 
What is provable about these two cardinals 
in ZF? The main result of this paper is that \ \ ZF 
$\vdash \forall A($\hbox{\carr$\Seq(A)$\carl}$
\neq$\hbox{\car${\cal{P}}(A)$\carl)  }
and we show that this is the best possible result. 
Furthermore, it is provable in ZF that if $B$  is an
infinite set, then  \car$
\FiniteB $\carl$<$\car${\cal{P}}(B)$\carl, 
even though the existence for some infinite set $B^{*}$  of  a function 
$f$ from $\em {fin}(B^{*})$ onto 
${\cal{P}}(B^{*})$ is consistent with ZF.
\end{abstract}
 
\vspace {1.0 cm}

%
%
%
%
{\large{\bf Section 0: }} {\em Introduction, Definitions and Basic Theorems} 

\vspace {0.3 cm}

{\em {Introduction:}}{   }In ZFC the cardinality of ordinal numbers 
plays an important 
role, since by AC 
each set has the cardinality
of some ordinal.\\
We use ``alephs'' for the cardinalities of ordinals. 
Thus in ZFC each cardinal number is 
an aleph. However this need not be the case
in ZF.\\
If we have a model $M$ of ZF in which the axiom of 
choice fails, then 
we have more cardinals 
in $M$ than in a model $V$ of ZFC,
even if we have fewer sets 
in $M$ than in $V$.
(This occurs when the choice-functions are not all in $M$).
This is because all ordinals are in $M$ and 
hence all alephs as well.\\
In this paper we are interested in the relation 
between three cardinals  arising in connection 
with a set $S$, namely,
\begin{eqnarray*}
1)&{\mbox{the cardinality of the power set of }}S\\
2)&{\mbox{the cardinality of the finite subsets of }}S\\
3)&{\mbox{the cardinality of the finite sequences 
         without repetition of }}S\\
\end{eqnarray*}
This section contains definitions and basic theorems provable in ZF.\\
In the next section we present two relative consistency proofs illustrating
 possible relations
between these cardinals.\\
The last two sections contain three results 
provable in ZF. The proofs of these are based on the same idea
originally from E. Specker, who used it to prove 
that the axiom of choice 
follows from the generalised continuum hypothesis [Sp1].
Assuming the existence of a function we derive a 
contradiction to Hartogs' Theorem.\hfill
\medskip

Because we do not use AC, our proofs are 
constructive. But we will see that 
sometimes arithmetic is powerful enough
for our constructions, making it an adequate substitute
for AC.\hfill\\
\smallskip 

{\em Cardinals:} A {\em cardinal number } $\cal {C}$ 
is the equvalence class of all sets which have the {\em{same size}}. 
(Two sets are said to have the same size {\em iff } 
there is a bijection between them.)\hfill
\medskip

{\em Alephs:} A cardinal number $\cal {C}$  is 
an {\em aleph} if it contains a well-ordered set. \\
We use calligraphic letters to denote cardinals and $\aleph$'s to denote 
the  alephs. \\
We denote the cardinality of the set $s$ by \car$s$\car.\hfill
\medskip

{\em Relations between cardinals:} We say that the cardinal number $\cal {C}$ 
is less than or equal to the cardinal number
$\cal {D}$ {\em iff } there are sets $c \in$ $\cal {C}$, $d \in$ 
$\cal {D}$ and a  \inj function 
from $c$ into $d$. \\
In this case we write $\cal {C}$ $\leq$ $\cal {D}.$
We write $\cal {C}$  $<$  $\cal {D}${  for  }
$\cal {C}$ $\leq$ $\cal {D}$  and  $\cal {C}$ $\neq$ 
$\cal {D}$. \hfill
\smallskip

If $c \in$ $\cal {C}$, $d \in$ 
$\cal {D}$ and we have a function 
from $d$ onto $c$, then we write 
$\cal {C} {\leq}^{*} \cal {D}$.\hfill
\smallskip

\bigskip

We  also need some well-known facts provable in ZF: \hfill
\bigskip

{\bf  Hartogs' Theorem:} 
Given a cardinal $\cal {C}$ there is a least 
aleph, $\aleph(\cal {C})$,  such that 
$\aleph(\cal {C}) \not\leq \cal {C}$ .\\
{\em {Proof:}} See [Je1] p.25 \sq \qed \\
\medskip

{\bf Cantor-Bernstein Theorem:} 
If \ls $\cal {C}$ and $\cal {D}$  {are cardinals with} 
$\cal {C}$ $\leq$ $\cal {D}$ and 
$\cal {D}$ $\leq$ $\cal {C}$, then $\cal {C}$ $=$ $\cal {D}$.\\
{\em {Proof:}} {See [Je1] p.23} \sq \qed \hfill
\pagebreak

{\bf Cantor Normal Form Theorem:}
{Any ordinal $\alpha$ can be written as} 
\[\alpha = \sum_{i=0}^{j} \omega^{\alpha_{i}}
\cdot k_{i}  \]   
{with }  $\alpha \geq \alpha_{0} > \alpha_{1} > \ldots > \alpha_{j} 
\geq 0 \ls, \ls \ls 1\leq k_{i} < \omega \ls, 
\ls \ls0 \leq j < \omega$.  \\
{\em {Proof:}} {See [Ba] p.57\halfsp\em{ff}} \sq \qed \\ 
\medskip

{\bf {Corollary 1: }}The Cantor Normal Form does not depend on AC.\\
{\em {Proof:}} {The proof of the Cantor Normal Form 
requires no infinite choices. \sq\qed}\\
\medskip

{\bf {Corollary 2: }}{If }\(\alpha = \sum\limits_{i=0}^{j} 
\omega^{\alpha_{i}} \cdot k_{i} \) 
{ is a Cantor Normal Form, then define 
$\stackrel{\hookleftarrow}{\alpha}$ by} 
\[\stackrel{\hookleftarrow}{\alpha} := \sum_{i=j}^{0} 
\omega^{\alpha_{i}}\cdot k_{i} = 
\omega^{\alpha_{0}} \cdot k_{0}.\]
{Then (in ZF) } \car$\alpha$\car$=$\car
$\stackrel{\hookleftarrow}{\alpha}$\car \\
{\em {Proof:}} {See [Ba] p.60} \sq \qed \\ 
\medskip

{\bf {Corollary 3: }} For any ordinal $\alpha$, 
ZF implies the existence of the following bijections.

\[\begin{array}{lll}
 F_{\Seq}^{\alpha}: &\ \alpha\longrightarrow
         \Seq (\alpha)&\ \ {\mbox{(=: finite sequences 
                                          of $\alpha$ without repetition)}}\\
 F_{\em{seq}}^{\alpha}\hfill: &\ \alpha\longrightarrow
   {\em{seq}}(\alpha)&\ \ {\mbox{(=: finite sequences 
                                          of $\alpha$\halfsp)}}\\
 F_{\em{fin}}^{\alpha}\hfill: &\ \alpha\longrightarrow
   {\em{fin}}(\alpha)&\ \ {\mbox{(=: finite subsets 
                                          of $\alpha$\halfsp)}}
\end{array}\]
{\em {Proof:}} Use the Cantor Normal Form Theorem, 
Corollary 2, order the finite subsets of $\alpha$ and then 
use the Cantor-Bernstein Theorem.\sq\qed\hfill 
\bigskip
\bigskip
\bigskip
 
%
%
%
%
{\large{\bf Section 1: }}{\em{ Consistency results}}

\vspace {0.3 cm}
 
In this section we work in the Mostowski permutation model to  
derive some relative consistency
results. The permutation models are models of  ZFA,  
set theory with atoms,
(see [Je2] p.44\halfsp{\em{ff}}\ls). \hfill\smallskip

The atoms $x \in A$ may also be considered to be sets which contain 
only themselves, this means: $x \in A \Rightarrow x=\{x\}$ 
(see [Sp2] p.197 or [La] p.2).\\
Thus the permutation models are 
models for ZF without the axiom of foundation. 
\hfill\smallskip
 
However, the Jech-Sochor Embedding Theorem 
(see [Je] p.208\halfsp{\em{ff}}\ls) implies
con\-sist\-ency results for ZF.\hfill
\medskip

In the permutation models we have a set of atoms $A$ and 
a group ${\cal G}$ of permutations of $A$.\hfill\smallskip

Let $\cal{F}$ be a normal filter on ${\cal G}$ (see [Je] p.199). 
We say that $x$ is {\em {symmetric}} if the group
${\mbox{sym}}_{\cal G}(x):=\{\pi \in {\cal G}:\ \pi (x)=x\ \}$ 
belongs to $\cal{F}$.\\
Let us further assume that ${\mbox{sym}}_{\cal G}(a) 
\in \cal{F}$ for every atom $a$, 
that is, that all atoms are symmetric 
(with respect to ${\cal G}$ 
and $\cal{F}$) and let $\cal{B}$ 
be the class of all hereditarily symmetric objects.\hfill\smallskip

The class ${\cal{B}}$ is both a permutation model and a transitive class:
all atoms are in ${\cal{B}}$
and $A\in{\cal{B}}$. Moreover, ${\cal{B}}$ is a 
transitive model of ZFA.\hfill\smallskip

Given a finite set $E \subset A$, let 
${\mbox{fix}}_{\cal G}(E):=\{\pi \in 
{\cal G}:\ \pi a=a $ for all $a\in E \}$
and let $\cal{F}$ be the filter on ${\cal G}$ generated by 
$\{{\mbox{fix}}_{\cal G}(E):\ E \subset A$ is finite$\}$.\\
$\cal{F}$ is a normal filter and $x$ is symmetric 
$\em{iff}$ there is a finite set of atoms $E_{x}$
such that $\pi(x)=x$ whenever $\pi \in {\cal G}$ and $\pi a=a$ 
for each $a \in E_{x}$. Such an $E_x$ is called a support for $x$.
\hfill

{\begin{tabbing}
Now \=the\= Mostowski model is constructed as follows: 
(see also [Je2] p.49\halfsp{\em{ff}}\ls)\\
        \>1) \>{The set of atoms} $A$ { is countable \\
        \>     \>(but need not be countable in the Mostowski model).}\\
        \>2) \>$R${ is an order-relation on }$A$.\\
        \>3) \>With respect to $R$, $A$ is a dense linear 
                 ordered set without endpoints.\\
        \>4) \>Let ${\mbox{Aut}}_R$ be the group of 
                    all permutations of $A$ such that\\
         \>     \>for all atoms $x,y\in A$, if $Rxy$ then 
                    $R\pi (x)\pi (y)$.\\
         \>5) \> Let ${\cal F}$ be generated by 
                       $\{\mbox{fix}(E):\  E\subset 
                       A \mbox{ is finite }\}$
\end{tabbing}
(Thus the atoms when ordered by $R$ are isomorphic to the 
rational numbers with the natural order).
We write $x<y$ instead of $Rxy$.\hfill\medskip

The subsets of $A$ (in the Mostowski model) are 
symmetric sets. 
Hence each subset of $A$ has
a finite support. \hfill
\smallskip

If $x \subseteq A$ (in the Mostowski 
model) and $x$ has non-empty support $E_{x}$, 
then an $a \in E_{x}$ may or may not
belongs to $x$.\hfill
\smallskip

\begin{tabbing}
{\bf {Fact: }}\=If $b \not\in x\cup E_{x}$ and there are two elements 
$a_{0},a_{1}\in E_{x}$ with $a_{0}<b<a_{1}$ such\\
\>that $\forall c(a_{0}<c<a_{1}\rightarrow c\not\in E_{x})$, then  
$\forall c(a_{0}<c<a_{1}\rightarrow c\not\in x)$.
\end{tabbing}
\begin{tabbing}
{\bf {Otherwise }}\=we construct a 
                          $\ \pi \in {\mbox{Aut}}_{R}$ such that  
                          $\pi a_{i}=a_{i}$ for 
		all $a_{i} \in E_{x}$ and $\pi c=b$.\\ 
                     \>Then $\pi (x) \neq x$, which is a contradiction.
\end{tabbing}
We can similarly show that if $a_{0}<b<a_{1}$ and $b \in x 
\setminus E_{x}$, then 
$\forall c(a_{0}<c<a_{1}\rightarrow c\in x)$. 
The cases when $\neg \exists a_{1}(b<a_{1}\ \wedge 
\ a_{1}\in E_{x})$
or $\neg \exists a_{0}(b>a_{0}\ \wedge \ a_{0}\in E_{x})$ are 
similar.\hfill
\smallskip

Hence, given a finite set $E \subset A$ (\car $E$\car=n$<\omega$), 
we can construct 
$2^{n} \cdot 2^{n+1}=2^{2n+1}$
subsets $x \subseteq A$ such that $E$ is a support of $x$.\hfill
\smallskip

Given a finite subset $E$ of $A$, consider the set ${\cal E}$ of
subsets of $A$ with support $E$. We use $R$ to order ${\cal E}$
as follows. Given $E_{1}=\{ a_{1},\ldots,a_{n}\}$ and 
$E_{2}=\{a_{1},\ldots,a_{n},\ldots,a_{n+k}\}$ 
with $a_{i}<a_{j}$ whenever $i<j$ and given $x\in {\cal E}$,
if $x$ is the $l^{th}$ subset with support $E_{1}$,
then $x$ is also the $l^{th}$ subset with support $E_{2}$.
\begin{tabbing}
Finally, we define the function {\em{Fin}}:  
\={\em{fin}}$(A)\ $\=$\longrightarrow \ {\cal{P}}(A)$ by\\
\>{    }$E$   \>$\longmapsto \ \ $\=\car $E$\car $^{th}$ subset of $A$ \\
\>                  \>	                            \> constructible with 
                                                           support $E$\\
\hfill\\
It is easy to see that {\em{Fin}} is onto.
\end{tabbing}
\pagebreak

If $E\subset A$ is finite, then use $R$ to order the
subsets of $E$ and use the corresponding lexicographic order on the
set of permutations of subsets of $E$. 
The set of permutations of subsets of
$E$ is isomorphic to $\Seq(E)$. In fact we 
can order $\Seq(E)$ for each finite $E\subset A$.\hfill
\smallskip

For each subset $x \subseteq A$ there is exactly one 
smallest support $E_{x}$ ($=:$supp$(x)$).\hfill\smallskip

If \hbox{\carr supp$(x)$\carl$=n$}, then put 
$\stackrel{\balk}{x}:=$\carr$\{y\subseteq A:
\ $supp$(y)=$supp$(x)\ \}$\carl$\leq 2^{2n+1}$ 
and for $l \leq \stackrel{\balk}{x}$ define as above
the $l^{th}$ element of
$\{y\subseteq A: \ $supp$(y)=$supp$(x)\}$. \\
We say that: ``$y \subseteq A$ is the $l^{th}$ 
subset of $A$ with support supp$(x)$''.\hfill\medskip

Now choose $24$ distinct elements $a_{0}, \ldots ,a_{23} \in A$ and 
define 
$A_{24}:=\{a_{0}, \ldots ,a_{23} \}$.\\ 
A simple calculation shows that 
\begin{center}
\hfill {{if }$n \geq 12$, then $2 \cdot 2^{2n+1} < n!$}\hfill$(*)$  
\end{center}
Take a finite subset $E$ of $A$ and let $y \subseteq A$ be the $l^{th}$ 
subset of $A$ with supp$(y)= E$.
Put $D:=${ supp}$(y)\Delta A_{24}$
{ (where $\Delta$ denotes 
symmetric difference)} and $d:=${\car}$D${\car.}\hfill
\smallskip

Define the function ${\em{Seq}}_{A}:\ {\cal{P}}(A)\longrightarrow 
\Seq (A)$ by
\[ {\em{Seq}}_{A}(y):= \left\{ \begin{array}{ll}
     {\mbox{ the }}l^{th}{\mbox{ permutation of supp}}(y) & 
     {\mbox{  if {\car}supp$(y)${\car}}}\geq 12,\\
     {\mbox{{ the }$(d!-l-1)^{th}${ permutation of supp}$(y)$}} &
     {\mbox{  otherwise.}}
      \end{array}
\right. \]

${\em{Seq}}_{A}$ is well defined because 
of $(*)$ and $d \geq13$.\hfill
\smallskip

It is easy to see that ${\em{Seq}}_{A}$ is \inj \backs . If there 
is a bijection between ${\cal{P}}(A)$ and 
$\Seq (A)$, then we find an 
$\omega$-sequence\ief in $A$ using an analogous construction. 
But this is a contradiction (see section 3).\\
\smallskip

Even more is true in the Mostowski model, (${\cal A}:=$\carl Atoms\carr),
\begin{center}
${\cal A} < {\em {fin}}({\cal A}) < {\cal{P}}({\cal A}) < \Seq ({\cal A}) < 
{\em{fin}}({\em{fin}}({\cal A})) < 
  {\em{seq}}({\cal A}) < {\cal{P}}({\cal{P}}({\cal A}))$.
\end{center}
(We omit the proof).\hfill
\bigskip

Our interest here is in the following result.
\begin{tabbing}
{\bf{Theorem 1:}} \=The following theories are equiconsistent:\\ 
\>(i)\ \ \ \ \=$\ {\mbox{ZF}}$\\ 
\>(ii)   \>$\ {\mbox{ZF}}+\exists {\cal A}({\cal{P}}
                     ({\cal A})< \Seq ({\cal A}))$\\ 
\>(iii)  \>$\ {\mbox{ZF}}+\exists {\cal A}({\cal{P}}
   ({\cal A}){\leq}^{*}{\em{fin}}({\cal A}))$
\end{tabbing}\medskip

{\em {Proof:}} It was shown above that in the Mostowski model 
there is a cardinal ${\cal A}$, namely the 
cardinality of the set of atoms, for which both (ii) and (iii) hold.\\
Unfortunately, the Mostowski model is only a model of ZFA. 
But it is well-known that 
Con(ZF)$\Rightarrow$Con(ZFC) and the Jech-Sochor 
Embedding Theorem provides 
a model of (ii) and (iii). \sq\qed \hfill\pagebreak

\begin{tabbing}
{\bf{Theorem 2:}} \=The following theories are equiconsistent:\\
\>(i)     \=$\ {\mbox{ZF}}$\\ 
\>(ii)   \>$\ {\mbox{ZF}}+\exists {\cal A}({\em {seq}}
                     ({\cal A})< {\em {fin}} ({\cal A}))$
\end{tabbing}
{\em {Proof:}}\hfill\\
By the Jech-Sochor Embedding Theorem it is enough
to construct a permutation model ${\cal B}$ in which 
there is a set $A$, such that:\hfill
\begin{tabbing}
    \=(a)  \=there is a \inj function from ${\em {seq}}(A)$ into 
                 ${\em {fin}}(A)$,\\
    \>(b)  \>there is no bijection between ${\em {seq}}(A)$ and
                 ${\em {fin}}(A)$.
\end{tabbing}\smallskip

We construct by induction on 
$n\in\omega$ the following:\hfill
\begin{tabbing}
    \=($\alpha$)  \=$A_0 :=\{\{\emptyset\}\};\ Sq_0 
                                :=(\{\emptyset\}):=\langle\ \rangle $;\\ 
    \>                       \>$G_0 :=$ the group of all permutations of 
                                $A_0$.
 \end{tabbing}
 Let $k_n$ be the number of elements of $G_n$, and ${\cal E}_n$ 
 be the set of sequences of $A_n$ of length less or equal than $n$ which
 does not belong to range$(Sq_n )$, then
 \begin{tabbing}
    \=($\beta$)     \=$A_{n+1} := A_n \dot\cup\{(n+1,\zeta,i):\ 
                                    \zeta\in {\cal E}_n \mbox{ and } 
                                    i<k_n + k_n\}$.\\
    \>($\delta$)     \>$Sq_{n+1}$ is a function from $A_{n+1}$ to ${\em
                                    {seq}}(A_n )$ defined as follows:
\end{tabbing}
   \[Sq_{n+1} (x)=\left\{\begin{array}{ll}
                                          Sq_n (x) &\mbox{if $x\in  
                                          A_n$,}\\
                                          \zeta &\mbox{if 
                                          $x=(n+1,\zeta,i)\in 
                                          A_{n+1}\setminus A_n$.}
                                          \end{array}\right. \]
\begin{tabbing}
    \=($\gamma$)   \=$G_{n+1}$ is the subgroup of the group of permutations
                                    of $A_{n+1}$ containing \\
     \>                       \>all permutations $h$ such that for some 
                                   $g_h \in G_n$ and $j_h <k_n + k_n$ we have 
\end{tabbing}                     
     \[ h(x)=\left\{\begin{array}{ll}
                                          g_h (x) &\mbox{if $x\in  
                                          A_n$,}\\
                                          (n+1,g_h (x), i{+}_n j_h) &\mbox{if 
                                          $x=(n+1,\zeta,i)\in 
                                          A_{n+1}\setminus A_n$.}
                                          \end{array}\right. \] 
\begin{tabbing}
\=($\gamma$)   \= \kill
 \> \>Where $g_h (\zeta)(m):=g_h (\zeta (m))$ 
                                   and ${+}_n$ is the addition modulo $k_n + 
                                   k_n$. 
\end{tabbing}

Let $A:=\bigcup\{A_n:\ n\in\omega\}$ and 
$Sq:=\bigcup\{Sq_n:\ n\in\omega\}$, then 
$Sq$ is a function from $A$ onto ${\em {seq}}(A)$.\hfill\smallskip

Further define for each natural number $n$ partial functions 
$f_n$ from $A$ to $A\cup \{\emptyset\}$ as follows. If $lg(x)$ denotes
the length of $Sq(x)$ and $n<lg(x)$, then $f_n (x):=Sq(x)(n)$, otherwise
let $f_n (x)=\emptyset$.\hfill\smallskip

Let Aut($A$) be the group of all permutations of $A$,
then ${\cal G}:=\{H\in \mbox{Aut}(A): \forall n\in\omega 
(H {\mid}_{A_n} \in G_n )\}$ 
is a group of permutations of $A$.
Let ${\cal F}$ be the normal filter on ${\cal G}$ generated by 
$\{\mbox{fix}(E):\ E\subset A \mbox{ is finite}\}$ and ${\cal B}$ be 
the class of all hereditarily symmetric objects.\\
Now $A\in {\cal B}$ and for each $n\in \omega$, supp$(f_n)=\emptyset$,
hence $f_n$ belongs to ${\cal B}$, too.\hfill\smallskip

Now define on $A$ a equivalence relation as follows,
$$x\sim y\ \ {\em{iff}}\ \ \forall n(f_n (x)=f_n (y)).$$
\pagebreak

{\bf{Facts:}}\hfill
\begin{enumerate}
\item Every equivalence class of $A$ is finite. (Because of each $A_n$ is
finite, hence each $k_n$).
\item${\em{seq}}(A)=\{\varsigma_x : x\in A\}$ where 
$\varsigma_x (n):=f_n (x)$, (if $f_n (x)\neq\emptyset)$.
\item For every finite subset $B$ of $A$, there are finite subsets $C,Y$ of
$A$ and a natural number $k>1$ such that $B\subseteq C$,  
$\ \forall x\in A\setminus C\ ($\car$\{H(x): H\in {\mbox{fix}}_{\cal G}(C)
\}$\carr$>k)$ and \carl$\{H[Y]: H\in {\mbox{fix}}_{\cal G}(C)\}$\carr$=k.$\\
(Choose $A_n \ (n\geq 1)$ such that $B\subseteq A_n$ and let $C:=A_n$.
Let $k:= k_n + k_n$ and $Y:=\{(n+1,\zeta,i)\in A_{n+1} :\ i\mbox{ is 
even}\}$. Then $Y$ has exactly two images under $\{h: \ h\in 
\mbox{fix}_{\cal G} (C)\}$ and for each $x\in A\setminus C$, 
\carl$\{h(x):\ h\in \mbox{fix}_{\cal G} (C)\}$\carr$\geq k_{n+1} +k_{n+1} 
>k$, because the $k_n$'s are strictly increasing.)
\end{enumerate}

Now the function\hfill
\[\begin{array}{rccc}
      \Psi\ :\ \ &{\em {seq}}(A)&\longrightarrow&{\em{fin}}(A)\\
                      &\varsigma&\longmapsto&\{x:\ {\varsigma}_x 
                      =\varsigma\}
\end{array}\]
is a \inj function in ${\cal B}$ from ${\em{seq}}(A)$ into 
${\em{fin}}(A)$ (by the facts 1 and 2).\\
Hence (a) holds in ${\cal B}$.\hfill\medskip

To prove (b), assume there is a \inj function $\Phi \in {\cal B}$ from
${\em{fin}}(A)$ into ${\em {seq}}(A)$.\hfill
\smallskip

Let $B$ be a support of $\Phi$ and let $C,Y,k$ be as in fact 3.
\hfill\smallskip

If the sequence $\Phi (Y)$ belongs to ${\em{seq}}(C)$, then for
some $H\in {\mbox{fix}}_{\cal G}(C),\ H[Y]\neq Y$, hence
$\Phi (H[Y])\neq \Phi (Y)$. But this contradicts that $H$ maps 
$\Phi$ to itself, (by definition of $C,Y$ and $H$). \hfill\smallskip

Otherwise there exists an $m\in \omega$ such that $x:=\Phi (Y)(m)$
does not belong to the set $C$.\\
Hence \carl$\{H(x):\ H\in {\mbox{fix}}_{\cal G}(C)\}$\carr$>k$ 
and \carl$\{H[Y]:\ H\in {\mbox{fix}}_{\cal G}(C)\}$\carr$=k$,
(by fact 3).\\
Every $H\in {\mbox{fix}}_{\cal G}(C)$ maps $\Phi$ to itself, hence
$\Phi (Y)$ to $\Phi (H[Y])$. So we have a mapping from a set with $k$
members onto  a set with more than $k$ members. \\
But this is a contradiction.\sq\qed
 
\vspace{3 cm}

%
%
%
%
{\large{\bf Section 2: }}  ${\mbox{ZF}} \vdash ($\car$\FiniteS$\car
		$<$\car$\PowerS$\car) for 
                              any infinite set $S$.

\vspace {0.5 cm}

{\bf{Theorem 3: }} ${\mbox{ZF}} \vdash 
{\em{fin}}{\cal{(C)}}<{\cal{P(C)}}$\hfill
\medskip

{\em {Proof:}} Take $S\in {\cal C}$. The natural map
from $\FiniteS$ into 
$\PowerS$ is a \inj function, hence 
\car$\FiniteS$\carl$\leq$\car$\PowerS$\carl
{ is always true.}\\
Assume that there is a bijective function 
\hbox{$B:\ \FiniteS \longrightarrow \PowerS$.} 
Then, given any  ordinal $\alpha$, we can construct 
an $\alpha$-sequence\ief in $\FiniteS$. But this contradicts Hartogs'
Theorem.\hfill
\pagebreak

First we construct an $\omega$-sequence\ief in 
$\FiniteS$ as follows:\hfill
\smallskip

$S \in \PowerS$ and, because $S$ is infinite, 
$S \not \in \FiniteS$.\\ 
But $B^{-1}(S) \in
\FiniteS$. So put $\so := B^{-1}(S)$ and 
$s_{n+1}:= B^{-1}(\sn)\ \ (n<\omega)$.\\
Then the set $\{\si:\ i<\omega\}$ is an infinite set of 
finite subsets of $S$ and the sequence 
$\langle \so, \se, \ldots ,\sn,\ldots {\rangle}_{\omega}$ is 
an $\omega$-sequence\ief in $\FiniteS$.\hfill
\medskip

If we have already constructed 
an $\alpha$-sequence\ief $\langle \so, \se, \ldots ,\sbe, \ldots 
{\rangle}_{\alpha}$ in $\FiniteS$ (with $\alpha \geq \omega $), 
then we define an equivalence relation on $S$ by
\begin{center}
$ x \sim y \ \ {\em{iff}}\ \ \forall\beta 
<\alpha(x\in\sbe\leftrightarrow y\in\sbe)$
\end{center}
Take $x \in S$ and suppose that $\mu < \alpha$. Define 
\begin{eqnarray*}
D_{x,\mu} & := & \bigcap_{\iota <\mu} \{\sio :\ x
	\in \sio\}\\
g(x)           & := & \{\mu <\alpha :\ x\in\smu\wedge 
	(\smu \cap D_{x,\mu}\neq D_{x,\mu})\}. 
\end{eqnarray*}

\begin{tabbing}
{\bf{Fact:}} \=Given $x,y \in S,\ g(x)=g(y) 
			\Leftrightarrow x \sim y.$\\
                      \>(In other words $x^{\sim}=y^{\sim}$ 
                                                      whenever $g(x)=g(y)$).\\
                      \>Hence  there is a bijection between 
                          $\{x^{\sim}: \ x \in S\}$ and $\{g(x):\ x \in S\}$.\\
                      \>Furthermore, $g(x) \in \final$.
\end{tabbing}
Since $\{g(x):\ x \in S\}\subseteq \final$, apply 
$F_{\em{fin}}^{\alpha}$ to obtain
$F_{\em{fin}}^{\alpha}[\{g(x):\ x \in S\}]\subseteq 
\alpha.$\hfill\smallskip

Let $\gamma$ be the order-type 
of $F_{\em{fin}}^{\alpha}[\{g(x):\ x \in S\}]$. Then 
$\gamma \leq \alpha$ and for each $g(x)$ we obtain
a number $\eta (g(x))<\gamma.$
\hfill
\smallskip

Each $\sio\ (\iota <\alpha)$ is the union of at most 
finitely many equivalence classes. Thus there is a \inj function 
\[\begin{array}{cccl}
h:\ &\alpha &\longrightarrow&\finga \\
      &\iota  &\longmapsto       &\{\xi:\ \eta (g(x))=\xi \wedge x\in \sio\}.
\end{array}\]
Since $F_{\em{fin}}^{\gamma}$ is a bijection between
$\finga$ and $\gamma$, 
$F_{\em{fin}}^{\gamma} \cirk h$ is a \inj function 
from $\alpha$ into $\gamma$ and because
$\gamma \leq \alpha$ we also have a \inj function 
from $\gamma$ into $\alpha$.\hfill 
\smallskip

The Cantor-Bernstein Theorem yields a bijection between
$\gamma$ and $\alpha$ and hence
a bijection $G$ from $\{\eta (g(x)): \ x\in S\}$ onto 
$\{\sio:\ \iota < \alpha \}$.\hfill
\smallskip

Now consider the function 
$\Gamma:= B\cirk G\cirk \eta \cirk g$ from $S$ into \PowerS:
$$\Gamma :\ S \stackrel{g}{\longrightarrow} \{g(x):\ x\in S\} 
     \stackrel{\eta}{\longrightarrow} \{\eta (g(x)): \ x\in S\}
     \stackrel{G}{\longrightarrow} \{\sio:\ \iota < \alpha \}
     \stackrel{B}{\longrightarrow} \PowerS $$

{\bf{Fact: }}$S_{\alpha}:=\ \{x\in S:\ x\not\in \Gamma (x)\} 
                        \not\in \{B(\sio):\ \iota <\alpha\}.$
\begin{tabbing}
{\bf{Otherwise }}\=Take $S_{\alpha}=B(\sbe)$ (for some $\beta < \alpha$). \\
                     \>We identify each $x^{\sim}$ with $g(x)$ using 
                     the bijection above.\\
                     \>Then there is a 
                         $g(x)$ such that $G\cirk \eta ((g(x))=\sbe$.\\
                      \>Now if $y \in x^{\sim}$ then $\Gamma 
                      (y)=S_{\alpha}$.\\
                       \>But $y \in S_{\alpha} 
\Leftrightarrow y \not \in \Gamma (y) 
\Leftrightarrow y \not\in S_{\alpha}$,
which is a contradiction.
\end{tabbing}

But $S_{\alpha} \subseteq S$ and 
$B^{-1}(S_{\alpha})=: \sa \in \FiniteS$ with 
$\sa \not\in \{\sio:\ \iota < \alpha \}$ and we have 
an $(\alpha +1)$-sequence\ief in $\FiniteS$, namely
$\langle \so, \se, \ldots 
,\sbe ,\ldots ,\sa {\rangle}_{\alpha +1}$. 
\medskip

We now see that for an infinite set $S$ there is no bijection
between $\FiniteS$ and $\PowerS$  and this completes
the proof.\sq\qed\hfill
\bigskip

We note the following facts.
\[\begin{array}{rl}
  {\mbox{Given a natural number $n$, }}&{\mbox{ZF}} 
		\vdash (n \times {\em{fin}}{\cal{(C)}}=
                                   		            {\cal{P(C)}}\rightarrow n=2^{k} 
                                                                    {\mbox{ for a $k<\omega $).}}\\
  {\mbox{Moreover, for each $k<\omega$  }}
                                      &{\mbox{Con(ZF)}}\Rightarrow 
		{\mbox{Con(ZF}}+\exists{\cal{C}}
                                            (2^{k}\times {\em{fin}}{\cal{(C)}}=
		{\cal{P(C)}})\\
                                  &{\mbox{(If $k=0$, then this is 
			obvious for finite cardinals.)}}
\end{array}\]
{\em{Sketch of the proof:}}\hfill\smallskip

For the consistency result, consider the permutation model with an infinite
set of atoms $A$ and the empty relation. Then the automorphism group is
the complete permutation group. It is not hard to see that any subset of 
$A$ in this model is either finite or has a finite complement. Take a natural number $k$
and consider (in this model) the set $k\times A$.
The cardinality of the set ${\cal{P}}(k\times A)$ is the same as that of 
the set $2^k \times\FiniteA$.\hfill\smallskip

To prove the other fact, assume that $n$ is a natural number which is not a power of $2$
and that for some infinite set $S$ there is a bijection $B$ between 
$n\times \FiniteS$ and ${\cal P}(S)$. Use the function $B$ to construct an 
$\omega$-sequence\ief in $\FiniteS$. Then, using Theorem 3, 
$\omega\leq \FiniteS < {\cal P}(S)$ and it is easy to see that $n\times \FiniteS
\leq \FiniteS \times \FiniteS\ =:{\FiniteS}^2.$\\
Then $\omega < {\cal P}(S) = n\times \FiniteS\leq{\FiniteS}^2$ 
contradicts the fact that 
if $\aleph_0\leq {\cal {P(C)}}$, then for any natural number $n$, 
${\cal {P(C)}}\not \leq {\em 
{fin}}({\cal C})^n$. (Here $\aleph_0$ 
denotes the cardinality of $\omega$). 
The proof of this fact is similar to the proof of Theorem 3.\sq\qed\hfill
\vspace{3 cm}
 
%
%
%
%
{\large{\bf Section 3: }}  $\Seq (S)$, {\em {seq}}$(S)$ and 
                        ${\cal{P}}(S)$ when $S$ is an arbitrary set.

\vspace {0.3 cm}

We show that  
${\mbox{ZF}} \vdash \Seq({\cal{C}}) \neq {\cal{P(C)}}
     {\mbox{  for every cardinal }}{\cal{C}}\geq 2$.
But we first need the following result. \hfill\medskip

{\bf {Lemma:}}  ${\mbox{ZF}} \vdash {\aleph}_{0} 
		\leq {\cal{P(C)}}\rightarrow 
                                         {\cal{P(C)}} \not \leq \Seq({\cal{C}})$.\hfill
\medskip

{\em {Proof:}}\\Take $S \in {\cal{C}}$. Then, because ${\aleph}_{0} \leq 
{\cal{P(C)}}$, we 
have a \inj function  
$\ f_{\omega}:\ \omega \longrightarrow \PowerS$.\\
Assume that there is a \inj function $J:\ \PowerS 
\longrightarrow \Seq (S)$.\\
Then $J\cirk f_{\omega}\backs :\omega \longrightarrow \Seq (S)$ is 
also \inj and we get an 
$\omega$-sequence\ief in $\Seq (S)$. Using 
this $\omega$-sequence\ief in $\Seq (S)$ we can easily
construct an $\omega$-sequence\ief in $S$.\hfill
\smallskip

If we already have constructed an $\alpha$-sequence\ief
$\langle \so,\se, \ldots ,\sbe,\ldots {\rangle}_{\alpha}$ 
($\alpha \geq \omega$) in $S$,
put $T:=\{\sio:\ \iota < \alpha\}$. This gives rise to
bijective functions, 
\[\begin{array}{crcl}
h_{0}:\ \ &T &\longrightarrow&\alpha\\
h_{1}:\ \ &\Seq (\alpha)&\longrightarrow&\Seq (T).
\end{array}\]
Let $J^{-1}$ be the inverse of $J$ 
and denote the inverse 
of $F_{\em{seq}}^{\alpha}$ by  $\mbox{inv}F_{\em{seq}}^{\alpha}$.
\hfill
\smallskip

Further define $$\Gamma := J^{-1}\cirk h_{1} \cirk 
		\mbox{inv}F_{\em{seq}}^{\alpha}
                           \cirk h_{0}$$
Note: dom$(\Gamma ) \subseteq T$ and range$(\Gamma) \subseteq \PowerS$ 
(because $J$ is \inj\backs).\hfill
\bigskip

{\bf{Fact: }}$S_{\alpha}:=\{x\in S:\ x\not\in\Gamma (x)\}
\not\in J^{-1}[\Seq (T)].$\hfill\smallskip

Assume not, then $x\in S$ such that 
$J({S}_{\alpha})=h_{1}\cirk 
    \mbox{inv}F_{\em{seq}}^{\alpha}\cirk h_{0}(x)$ 
 yields a contradiction.\hfill\medskip
    
Because $J({S}_{\alpha}) \not\in \Seq (T)$, 
the sequence $J({S}_{\alpha})$ has a first element which is not in $T$, say  $\sa$.
Finally, the sequence $\langle \so,\se, \ldots ,\sa 
{\rangle}_{\alpha +1}$ is 
an ($\alpha +1$)-sequence\ief in $S$.
\hfill\smallskip

So the existence of a \inj function $J:\ \PowerS 
\longrightarrow \Seq (S)$ contradicts Hartogs' Theorem.\sq\qed 
\hfill\bigskip

{\bf {Theorem 4:}}\hfill If ${\cal C}\geq 2$ is any cardinal, then
${\mbox{ZF}}\vdash (\Seq {\cal{(C)}} \neq {\cal{P(C)}})$\hfill\medskip

{\em {Proof:}}\hfill\\
By the Lemma it is enough to prove that if ${\cal{C}}\geq 2$, then 
                                          $\Seq ({\cal{C}})={\cal{P(C)}}\Rightarrow  
                                           {\aleph}_{0}\leq{\cal{C}}$.\hfill
\smallskip

For finite cardinals ${\cal{C}}\geq 2$ the statement is obvious.
So let $S\in{\cal{C}}$ be an infinite set and assume that there 
is a bijective function
$$B: \Seq (S) \longrightarrow {\cal{P}}(S).$$
We use this function to construct an $\omega$-sequence\ief in $S$.\hfill
\smallskip

Let ${n}^{\star}\ (n<\omega)$ be the cardinality of $\Seq (n).$\\
Then ${0}^{\star}=1;\ {1}^{\star}=2;\ {2}^{\star}=5;
	\ldots \ {16}^{\star}=
                 56,874,039,553,217;\ldots$
                \ \hbox{(see [Sl], No. 589),}\\
and, in general   \[{n}^{\star}=\sum_{i=0}^{n} \frac{n!}{i!}\]

We begin by choosing four distinct elements 
of $S,\ S_{4}:=\{s_{0},s_{1},s_{2},s_{3}\}$ and
use these elements to construct 
a $4$-sequence\ief $\langle s_{0},s_{1},s_{2},s_{3}{\rangle}_{4}$
in $S$. This sequence will give us 
an order on the set $\Seq ({S}_{4})$ (e.g. we order 
$\Seq ({S}_{4})$ by length and lexicographically).\hfill
\smallskip

If we have already constructed an $n$-sequence\ief 
$\langle s_{0},s_{1},\ldots,s_{n-1}{\rangle}_{n}$ in $S$ 
($n\geq 4$), then put
$S_{n}:=\{s_{i}:\ i<n\}$.
Then $B[\Seq ({S}_{n})]\subseteq {\cal{P}}(S)$ has cardinality 
${n}^{\star}$.\hfill
\smallskip

We now define an equivalence relation on $S$ by
$$x\sim y \ \ {\em{iff}}\ \ \forall q
\in \Seq ({S}_{n})(x \in B(q) \leftrightarrow y \in B(q)).$$
It is easy to see that for each $q \in \Seq ({S}_{n})$
\begin{eqnarray}
B(q){\mbox{ is the {\em{disjoint union}} of }}
{\mbox{{\em{less}} than ${n}^{\star}$ equivalence classes.}}
\end{eqnarray}
Take the above order on $\Seq ({S}_{n})$. This induces 
an order on the set of equivalence
classes eq$:=\{x^{\sim}:\ x\in S\}$ and also an order 
on ${\cal{P}}$(eq).\hfill
\pagebreak

If there is a first $r \in {\cal{P}}({\mbox{eq}})$ such that 
$r \not\in B[\Seq ({S}_{n})]$,
then $q_{r}\backs :={B}^{-1}(r)$ is a ``new'' sequence in S. 
This is $q_{r}\not\in\Seq ({S}_{n})$
and we choose the first element $\sn$ of $q_{r}$ which is 
not in $S_{n}$.\\
Hence, the sequence $\langle \so,\se,\ldots ,\sn {\rangle}_{n+1}$ is 
now an $(n+1)$-sequence\ief in $S$.
\hfill
\smallskip

If there is an $\si \in {S}_{n}$ such that 
$\{\si\}\not\in B[\Seq ({S}_{n})]$, then use 
$B(\{\si\})$ to construct an $(n+1)$-sequence\ief in $S$.\hfill
\smallskip

Otherwise our construction stops at $S_{n}$ and we
write stop($S_{n}$).\hfill\medskip

Our construction only stops if
\[\begin{array}{ll} 
{\mbox{for each }}\si \in S_{n}:&\ \{\si\}\in {\mbox{eq  and }}\\
{\mbox{for each }}r\in {\cal{P}}{\mbox{(eq)}}&
\ {\mbox{there is a }}q_{r} \in \Seq ({S}_{n})
{\mbox{ such that }}B({q}_{r})=r.
\end{array}\]
If $\kappa\ (\kappa<\omega)$ is the cardinality of eq, 
then $2^{\kappa}$ is the 
cardiniality of ${\cal{P}}$(eq) and because of
$(1)$ we have 
${\mbox{stop}}({S}_{n})\ \rightarrow \ 2^{\kappa}={n}^{\star}.$\\ 
It is known that ${{0}^{\star}}=1={{2}^{0}};\ {{1}^{\star}}=2={{2}^{1}};\ 
                     {{3}^{\star}}=16={{2}^{4}}$ and 
                     ${{n}^{\star}}$ is a power of $2$ for some $n>3$,
then $n$ has to be bigger than ${10}^{8}$.\hfill\smallskip

If there are only finitely many  $k,n<\omega$ such that
${2}^{k}={n}^{\star}$, then
there is a least ${n}_{0}$ such that ${2}^{k}={{n}_{0}}^{\star}$ 
and $\forall n>{n}_{0}(\neg {\mbox{stop}}({S}_{n})).$\hfill
\bigskip

Refining our construction removes the need for this 
strong arithmetic condition.
\hfill\medskip

Assume stop($S_{n}$).\\
If $x\not\in S_{n}$ then let $S_{n+1}^{x}:= {S}_{n}\dot\cup\{x\}$ and 
$S_{n+k}^{x}:= {S_{n+1}^{x}}\dot\cup\{Y\}$ with $Y$ of cardinality $k-1$.
Because  ($n$ is even)$\Leftrightarrow ({n}^{\star}$ is odd){    }
and stop($S_{n}$),
we cannot have  stop($S_{n+1}^{x}$)
for any $x \not\in {S}_{n}$.\hfill\medskip

Now we recommence our construction 
with the set $S_{n+1}^{x}$ and construct an 
$(n+k)$-sequence\ief $\langle \so,\se,\ldots,{s}_{n+k-1}
{\rangle}_{n+k}\ (k\geq 2)$ in $S$.\\
Define $S_{n+k}^{x}:=\{\si:\ i<(n+k)\}$.
If this construction also stops at the
$(n+{\mbox{\footnotesize{stop}}})^{th}$ stage at the set
$S_{n+{\mbox{\scriptsize{stop}}}}^{x} \ \ (
{\mbox{\footnotesize{stop}}}\geq 2)$,
then we write  ${S}^{x}$  instead of  
${S}_{n+{\mbox{\scriptsize{stop}}}}^{x}$.\hfill\medskip

If there is an $x \in S$ such that $S^{x}$ is infinite, then 
our construction does not stop when
we recommence with ${S}_{n+1}^{x}$ and
we can construct an $\omega$-sequence\ief in $S$. But 
this contradicts our Lemma.
\hfill\medskip

So there cannot be such an $x$ and each $x\in S$ is in 
exactly one {\em{finite}} set $S^{x}$.
If for each $x\in S,\ S^{x}$ is the union of some elements of eq, 
then $S$ must be finite, because eq is finite. But this 
contradicts our assumption that $S$ is infinite.\hfill

\begin{center}
A subset of $S$ is called {\em {good}} if it cannot be written as the
union of elements of eq.\hfill
\end{center}

Consider the set $T_{\mbox{\scriptsize{min}}}:=\{x:\ {S}^{x} 
{\mbox{ is good and of least cardinality}}\}$
and let ${{m}_{\mbox{\tiny{T}}}}$ be the cardinality of ${S}^{x}$ for 
some $x$ in 
$T_{\mbox{\scriptsize{min}}}$. 
Further for $x\in T_{\mbox{\scriptsize{min}}}$ let $x_= :=
\{y:\ S^y = S^x\}$ (this elements of $S^x$ we cannot distinguish) 
and $m_=$ denote the least cardinality of the sets $x_=$.\hfill
\smallskip

If $T_{\mbox{\scriptsize{min}}}$ is good, use
$B^{-1}(T_{\mbox{\scriptsize{min}}})$ to construct
an $(n+1)$-sequence\ief in $S$.\hfill\smallskip

Otherwise take $x\in T_{\mbox{\scriptsize{min}}}.$
Because $S^x$ is good
$${B}^{-1}({S}^{x})\not\in \Seq ({S}_{n}).$$
Thus there is a first $y$ in ${B}^{-1}({S}^{x})$ which 
is not in $S_{n}$. It is easy to see that $S^{y}\subseteq S^{x}$
and if $S^{y}\neq S^{x}$ then  
$S^{y}$ is not good (because of $x \in T_{\mbox{\scriptsize{min}}}$).\\
But then $B^{-1}({S}^{x}\setminus {S}^{y})\not\in \Seq ({S}^{y})$ and 
we  may proceed.
\hfill\medskip

So for each $x\in T_{\mbox{\scriptsize{min}}}$
construct an ${{m}_{\mbox{\tiny{T}}}}$-sequence\ief 
${\mbox{SEQ}}^{x}$
in $S$ such that 
$$S^x =S^y\ \Longrightarrow 
\ {\mbox{SEQ}}^{x}={\mbox{SEQ}}^{y}.$$
For $i<{{m}_{\mbox{\tiny{T}}}}$ define
$$Q_{i}:=\{s \in S:\ {\mbox{$s$ is 
		the $i^{th}$ element in }}{\mbox{SEQ}}^{x}
                                   {\mbox{ for some }}x\in S\}$$
Assume there is some $j<{{m}_{\mbox{\tiny{T}}}}$ such 
that $Q_{j}$ is good. Then 
$B^{-1}({Q}_{j})\not\in\Seq ({S}_{n})$. But 
$B^{-1}({Q}_{j})\not\in\Seq (S)$ and we 
get an $(n+1)$-sequence\ief in $S$.\hfill
\medskip

It remains to justify our assumption.\hfill\smallskip

Note that if for some $i\neq j$, $z\in Q_i \cap Q_j$,  then $S^z$ cannot 
be good and $z\not\in\cup\{x_= : x\in T_{\mbox{\scriptsize{min}}}\}$.
Furthermore for each $x \in T_{\mbox{\scriptsize{min}}}$ there is exactly
one $i_x$ such that $x\in Q_{i_x}$ and if $z,y \in x_= ,\ z\neq y$, then
$i_x \neq i_y$. If there are no good $Q_i$'s, $m_=$ cannot exceed
 $\kappa$, (the cardinality of eq). But by the following this is a 
contradiction:\hfill\medskip

An easy calculation modulo $2^{r}\ (r\leq 4)$ shows that for each $n$,
if $2^{r}|{n}^{\star}$, then $2^{r}|{(n+{2^{r}})}^{\star}$ and
${2^{r}}\goback\not |\halfsp{(n+t)}^{\star}$ if $0<t<{2^{r}}$.
\hfill\medskip

Assume there is a smallest $k\ (k\geq 4)$ such that $2^{k+1}|{n}^{\star}\  
                            {\mbox{and}}\ 2^{k+1}|{(n+t)}^{\star}\ 
                            {\mbox{for some }}t \\{\mbox{with }}
                                     0<t<{2^{k+1}}.$\hfill
\smallskip

Then, because $2^{k}|2^{k+1}$, we have 
$2^{k}|{n}^{\star} \ {\mbox{and}}\ 2^{k}|{(n+t)}^{\star}$. 
Since $k$ is by definition the smallest 
such number, we know that $t$ must be $2^{k}$.\hfill
\[\begin{array}{lrcrcrr}
{(n+{2^{k}})}^{\star}={\sum\limits_{i=0}^{n+{2^{k}}} 
\frac{(n+{2^{k}})!}{i!}}=&
 \ \ 1\cdot 2\cdot&\ldots&\cdot
{2^{k}}\cdot &\hspace{-3mm}
   ({2^{k}}+1)\cdot\ \ldots&\cdot({2^{k}}+n)&
{\mbox{\scriptsize{(1)}}}\\
 &+  \hfill 2\cdot&\ldots&\cdot {2^{k}}\cdot &\ldots                              
 &\cdot({2^{k}}+n)&{\mbox{\scriptsize{(2)}}}\\
 &&\ddots&&&&\vdots \ \\
 &+\hfill &&{2^{k}}\cdot&\ldots &\cdot({2^{k}}+n)&
 {\mbox{\scriptsize{(${2^{k}}$)}}}\\
 &&&&\ddots &&\vdots\ \\
 &+\hfill &&&&({2^{k}}+n)&{\mbox{\scriptsize{(${2^{k}}+n$)}}}\\
 &+\hfill &&&&1\ &{\mbox{\scriptsize{(${2^{k}}+n+1$)}}}
\end{array}\]

It is easy to see that $2^{k+1}$ divides 
lines {\scriptsize{$(1)-({2^{k}})$}}
since $k \geq 2{\mbox{ and }}n\geq2$.\hfill\smallskip

If we calculate the products of lines 
{\scriptsize{$({2^{k}}+1)-({2^{k}}+n+1)$}},
then we only have to consider 
sums which are not obiviously divisible by $2^{k+1}$. So, for a
suitable natural number ${\varepsilon}$ we have\hfill
\begin{eqnarray}
{(n+{2^{k}})}^{\star}={2^{k}}\cdot (\sum_{j=0}^{n-1}
   \sum_{i>j}^{n} \frac{n!}{i\cdot j!})+{n}^{\star} + {2^{k+1}}\backs
   \cdot \varepsilon .
\end{eqnarray}
{We know that }$2^{k+1}|{n}^{\star}${  with  }$n\geq 3,\ k\geq 4 $.
And because $n^{\star}$ is even $n{\mbox{ has to be odd.}}$\\
Then for $j=\ n-1,\ n-2$ or 
$n-3 \ \sum\limits_{i>j}^{n} \frac{n!}{i\cdot j!}$ is odd. Moreover, 
if $0 \leq j\leq (n-4)$, then
$\sum\limits_{i>j}^{n} \frac{n!}{i\cdot j!}$
is even. ${\mbox{So   }} \sum\limits_{j=0}^{n-1} \sum\limits_{i>j}^{n} 
\frac{n!}{i\cdot j!}{\mbox{ is odd. Hence  }}
{2^{k+1}}\goback\not |\halfsp {(n+{2^{k}})}^{\star},
{\mbox{  (by $(2)$ and $2^{k+1}|{n}^{\star}$).}}$
\bigskip
\medskip

We return to the proof.\hfill\smallskip

We know that if $2^k = n^{\star}$ and $(n+t)^{\star}$ is a power of $2$,
then $2^k$ divides $t$.\hfill $(**)$\hfill\smallskip

Take $x\in T_{\mbox{\scriptsize{min}}}$ such 
that \carl$x_=$\carr$= m_=$. If $y\in S^x$, then \hfill
\begin{tabbing}
\ \ \ \ \ \ \ \ \ \ \ \ \=(i)\ \ \=\carl$S^y$\carr$=n+t_y$ 
                   with $2^{\kappa}$ divides $t_y$,\\
                  \> (ii) \>either $y\in x_=$ or $S^y$ is not good.
\end{tabbing}

This is because $2^{\kappa} = n^{\star}$ and $(**)$.\hfill\smallskip

Hence (for a suitable natural number $\varepsilon$) 
$m_{\mbox{\tiny{T}}} =$\carl
$S^x$\carr$=n+2^{\kappa}\cdot\varepsilon +m_=$ 
(by (ii)), and $2^{\kappa}$ divides $m_=$ (by (i)). \\
But this implies that $m_=$ must be larger 
than $\kappa$, which justifies our
assumption.\sq\qed\hfill\\            
\bigskip

The statement obtained when $\Seq$ is replaced by {\em {seq}} is much easier
to prove:\hfill\bigskip
                       
{\bf {Theorem 5: }} ZF $\vdash \ {\em {seq}}
({\cal C})\neq {\cal P}({\cal C})$ for all cardinals such that $\emptyset
\not\in {\cal C}.$\hfill\smallskip

{\em {Proof: }} Take $S\in {\cal C}$. First note the fact that 
if $\aleph_0 \leq {\cal C}$, then ${\em {seq}}({\cal C}) \not \geq
{\cal P}({\cal C})$.\hfill\smallskip

(The proof is the same as the proof of the Lemma, except that
we can skip the first lines of the proof of the Lemma).\hfill\smallskip

Assume there is a bijection $B$ from ${\em {seq}}(S)$ onto ${\cal P}(S)$.
Choose an $s_0 \in S$, and define a \inj function $f_{s_0}$ from
$\omega$ into ${\cal P}(S)$ by $i\mapsto \xi_i:=B(\langle s_0,s_0,
\ldots ,s_0 \rangle)\ (i$-times). Use the $\xi_i$'s to construct pairwise
disjoint subsets $c_i \subseteq S \ (i<\omega)$.\\
Given an $n$-sequence\ief $\langle s_0,s_1,\ldots,s_{n-1}{\rangle}_n$
in $S$, let $S_n :=\{s_i : \ i<n \}$ and the natural order on $S_n$ 
induce a well-ordering on the set ${\em {seq}}( S_n )$ with order type 
$\omega$. Then there is a bijection $h:\ \omega \longrightarrow {\em 
{seq}}( S_n )$.\\
The function $\Gamma := B\cirk h$ is a \inj function from $\omega$ into 
${\cal P}(S)$ and $t:=\dot\cup\{c_i : c_i \subseteq \Gamma (i)\}\not\in
\{\Gamma (k): k<\omega\}$.\\
Hence $B^{-1} (t)$ is a sequence in $S$ which does not belongs to $S_n$. 
Choose $s_n \in S$ to be the first element of $B^{-1} (t)$ not in $S_n$. Then
$\langle s_0,s_1,\ldots,s_n {\rangle}_{n+1}$ is an $(n+1)$-sequence\ief
in the set $S$.\hfill\smallskip

We thus construct an $\omega$-sequence\ief in $S$, contradicting 
the previous fact.\sq\qed\hfill\bigskip

\pagebreak
\bf {References}
\begin{description}
         \item   [{[Ba]}] H.\halfsp Bachmann, Transfinite Zahlen, 
                                            Springer-Verlag, Berlin, 1967
         \item [{[Je1]}] Th.\halfsp Jech, Set Theory, Academic Press, 
                                            New York, 1978
         \item [{[Je2]}] Th.\halfsp Jech, The Axiom of Choice, 
		       North-Holland Publ.\halfsp Co., 
                                   Amsterdam, 1973
         \item   [{[La]}] H.\halfsp L\"auchli, Auswahlaxiom in der Algebra, 
		       Comment.\halfsp Math.\halfsp Helv., 
                                  vol.37, 1962, pp.1-18
         \item   [{[Sl]}] N.J.A.\halfsp Sloane, A Handbook of Integer Sequences, 
		      Academic Press, 
                                New York, 1973
         \item[{[Sp1]}] E.\halfsp Specker, Verallgemeinerte 
		     Kontinuumshypothese und 
                                    Aus\-wahl\-axiom, Archiv der Mathematik 5, 
		     1954, pp.332-337
         \item[{[Sp2]}] E.\halfsp Specker, Zur Axiomatik der Mengenlehre, 
		    Zeitschr.\halfsp f.\halfsp math.\halfsp 
                                   Lo\-gik und Grundl.\halfsp der 
                                   Math. 3, 1957, pp.173-210
\end{description}
\end{document}